\documentclass[12pt,reqno]{article}
\usepackage{amsfonts}
\usepackage{amsthm}
\usepackage[mathscr]{eucal}
\usepackage{amsmath}
\usepackage{amssymb}
\usepackage{amscd}
\makeatletter
\def\newtheorem#1{\@ifnextchar[{\@othm{#1}}{\@nthm{#1}}}
\def\@nthm#1#2{%
\@ifnextchar[{\@xnthm{#1}{#2}}{\@ynthm{#1}{#2}}}
\def\@xnthm#1#2[#3]{\expandafter\@ifdefinable\csname #1\endcsname
{\@definecounter{#1}\@newctr{#1}[#3]%
\expandafter\xdef\csname the#1\endcsname{\expandafter\noexpand
  \csname the#3\endcsname \@thmcountersep \@thmcounter{#1}}%
\global\@namedef{#1}{%
  \@thm{#1}{#2}}\global\@namedef{end#1}{\@endtheorem}}}
\def\@ynthm#1#2{\expandafter\@ifdefinable\csname #1\endcsname
{\@definecounter{#1}%
\expandafter\xdef\csname the#1\endcsname{\@thmcounter{#1}}%
\global\@namedef{#1}{%
  \@thm{#1}{#2}}\global\@namedef{end#1}{\@endtheorem}}}
\def\@othm#1[#2]#3{%
  \@ifundefined{c@#2}{\@nocounterr{#2}}%
  {\expandafter\@ifdefinable\csname #1\endcsname
  {\global\@namedef{the#1}{\@nameuse{the#2}}%
\global\@namedef{#1}{\@thm{#2}{#3}}%
\global\@namedef{end#1}{\@endtheorem}}}}
\def\@thm#1#2{\refstepcounter
    {#1}\@ifnextchar[{\@ythm{#1}{#2}}{\@xthm{#1}{#2}}}
\def\@xthm#1#2{\@begintheorem{#2}{\csname the#1\endcsname}\ignorespaces}
\def\@ythm#1#2[#3]{\@opargbegintheorem{#2}{\csname
       the#1\endcsname}{#3}\ignorespaces}
\def\@thmcounter#1{\noexpand\arabic{#1}}
\def\@thmcountersep{.}
\def\@begintheorem#1#2{\trivlist
   \item[\hskip \labelsep{\bfseries #2\ #1.}]\itshape}
\def\@opargbegintheorem#1#2#3{\trivlist
      \item[\hskip \labelsep{\bfseries #2\ #1\ [#3].}]\itshape}
\def\@endtheorem{\endtrivlist}
\makeatother
\theoremstyle{plain}
\newtheorem{lemma}{Lemma}[section]
\newtheorem{theorem}[lemma]{Theorem}
\newtheorem{proposition}[lemma]{Proposition}

\theoremstyle{definition}

\def\d{$\displaystyle}
\def\be{\begin{equation}}
\def\heavycdot{\raisebox{2pt}{\tiny\kern.5em$\bullet$}}

\numberwithin{equation}{section}
\def\begeq{\stepcounter{lemma}\begin{equation}}

\date{}

\begin{document}

\title{Cyclic Homology of $DG$ Coalgebras and a K\"{u}nneth Formula}
\author{Masoud Khalkhali\thanks{\; Author's research is partially supported
by NSERC of Canada.  Presented in a conference in honor of Professors
Peter Fillmore and Heydar Radjavi, Halifax, Nova Scotia, June 1998.}}

\maketitle

%%%%%%%%%%%%%% abstract

%\begin{abstract}
%%%%%%% put the text of the abstract here
%\end{abstract}

In this paper we extend the cyclic homology functor, and in particular
the periodic cyclic homology, to the category of $DG$ (= differential
graded) coalgebras.  We are partly motivated by the question of products
and coproducts in periodic cyclic homology of algebras.  As an
application, we will show how one can start from the classical shuffle
map in homological algebra and algebraic topology, interpreted in [HMS]
as a morphism of $DG$ coalgebras, and build a theory of products and
coproducts in (periodic, negative, etc.) cyclic homology.  Along these
lines we also recover a K\"{u}nneth isomorphism
\[ \widehat{S}:CC_\ast (A_1)\widehat{\otimes}CC_\ast
(A_2)\overset{\sim}{\longrightarrow}CC_\ast (A_1\otimes A_2)\; , \]
relating the periodic cyclic complexes of unital $(DG)$ algebras over a
field of characteristic zero.  Here $CC_\ast$ denotes the periodic
cyclic complex and $\widehat{\otimes}$ denotes the topologically
completed tensor product of complexes.  We note that a homotopy
equivalent result, with a different choice for $CC_\ast$ and
$\widehat{S}$, can be
found in the works of Cuntz and Quillen [CQ$_2$], Puschnigg [P] and more
recently Bauval [B].  Our method is based on ideas of Cuntz and Quillen
in cyclic homology (cf. [CQ$_1$] and reference therein) and in a sense
is dual to them.

An interesting problem suggests itself.  It would be interesting to see
how one can use a similar approach to define the cyclic cohomology of
$(DG)$ Hopf algebras and compare it with the definition of periodic
cyclic cohomology of Hopf algebras as recently defined by Connes and
Moscovici [CM].

I am much indebted to Joachim Cuntz for his comments on a first draft of
this paper and ensuing discussions which greatly improved my understanding
of the K\"{u}nneth formula and the importance of the topological tensor
product.  I would also like to thank A. Bauval for communicating her
results and for very informative discussions.  I am also grateful to J.
Williams for typing the manuscript.

\section{Bar and cobar constructions and the shuffle map}
The present form of the bar and cobar constructions first appeared in
a paper by Husemoller, Moore and Stasheff [HMS] under the name of
{\em algebraic classifying space} and {\em loop space constructions}.
The reason for these names is that these constructions provide models
for singular chains on the classifying space of groups and singular cochains
on the (based) loop space of simply connected spaces.  In the same paper the
{\em shuffle map} is introduced as a morphism of $DG$ coalgebras.  This is
crucial for applications to periodic cyclic homology.  The bar construction,
denoted here by $B$, is a functor from the category of $DG$ algebras to the
category of $DG$ coalgebras.  It has a left adjoint, namely the cobar
construction, denoted here by $B^c$, from the category of $DG$ coalgebras to
$DG$ algebras.  Let us recall their definitions.

Let \d A = \oplus_{i\geq 0}A_i$ be a positively graded $DG$ algebra
over a ground ring $k$.  We assume its differential $d$ has degree $-1$.
Let $A[1]$ denote the suspension of $A$ defined by $A[1]_n=A_{n-1}$.
The {\em bar construction} of $A$, denoted by $BA$, is a $DG$ coalgebra whose
underlying graded coalgebra is $BA=T^cA[1]$, the cofree coaugmented
counital coalgebra generated by the graded vector space $A[1]$.  We have
$(BA)_0=k$, $(BA)_1=A_0$ and in general
\[ (BA)_n = \bigoplus^n_{s=0}\;\bigoplus_{i_1+\cdots +i_r+r=s}\;
A_{i_1}\otimes\cdots\otimes A_{i_r}\; . \]
The differential of $BA$ is $d+b^\prime$, where $d$ and $b^\prime$ are
defined by
\[ d(a_1,\dots ,a_n) = \sum^n_{i=1}(-1)^{\varepsilon_i}(a_1,\dots
,da_i,\dots ,a_n)\]
and
\[ b^\prime (a_1,\dots ,a_n) =
\sum^{n-1}_{i=1}(-1)^{\varepsilon_i+|a_i|}(a_1,\dots ,a_ia_{i+1},\dots
,a_n)\; , \]
where
\[ \varepsilon_i = \sum^{i-1}_{j=1}|a_j|+i-1\; . \]
Computation shows that $d^2 = b^{\prime\; ^2} = db^\prime +b^\prime d=0$,
so that $d+b^\prime$ is a differential of degree -1.  It is moreover a
graded coderivation of $BA$.

Morphisms of $DG$ coalgebras into $BA$ are defined via {\em twisting
cochains} and vice versa.  More precisely, let $C$ be a $DG$ coalgebra.  A
degree -1 map $\theta :C\to A$ is called a twisting cochain if
\[ \delta\theta +\theta^2 = 0\; , \]
where, $\delta\theta = [d,\theta ]$ in the Hom complex differential of
$\theta$ and $\theta^2= m\circ (\theta\otimes\theta)\Delta$, where
$\Delta$ is the coproduct of $C$ and $m$ is the multiplication of $A$.
The {\em universal twisting cochain} is the map $\theta :BA\to A$ defined by
$\theta (a) = a$ for $a\in A$ and zero otherwise.  There is a 1-1
correspondence between morphism of $DG$ coalgebras $\widehat{\theta}:C\to
BA$ of degree zero and twisting cochains $\theta :C\to A$.  In this
correspondence, $\theta$ is simply the corestriction of
$\widehat{\theta}$ to $A$.  Conversely, given $\theta$, we have
\begin{equation*}
\widehat{\theta} = \sum^\infty_{n=1}\widehat{\theta}_n\tag{1}
\end{equation*}
where
\[\widehat{\theta}_n = \theta^{\otimes n}\circ\Delta^{(n)}\; , \]
and $\Delta^{(n)}$ is the $n-th$ iteration of the coproduct.

Next, we discuss the {\em cobar construction}.  Let \d C =
\oplus_{i\geq 0}C_i$ be a positively graded $DG$ coalgebra with a
differential $d$ of degree -1.  We assume $C_0 = k$.  Let $C[-1]$ be the
graded space defined by $C[-1]_i=0$ if $i<0$, and $C[-1]_i=C_{i+1}$ for
$i\geq 0$.  The cobar construction of $C$, denoted by $B^cC$, is the
$DG$ algebra whose underlying graded algebra is $B^cC=TC[-1]$, the augmented
unital free algebra generated by $C[-1]$.  We have \d (B^cC)_0 =
\oplus_{n\geq 0}C^{\otimes n}_1$, and
\[ (B^cC)_n = \bigoplus^n_{r=1}\;\bigoplus_{i_1+i_2+\cdots +i_r=n+r}\;
C_{i_1}\otimes\cdots\otimes C_{i_r}\; . \]
The differential of $B^cC$ is $d+b^\prime$, where $d$ and $b^\prime$ are
defined by
\[ d(a_1,\dots , a_n) = \sum (-1)^{\varepsilon_i}(a_1,\dots , da_i,\dots
, a_n) \]
and
\[ b^\prime = \sum
1\otimes\cdots\otimes\tilde{\Delta}\otimes\cdots\otimes 1\; , \]
where $\varepsilon_i$ is the same as before and
$\tilde{\Delta}:C[-1]\longrightarrow C[-1]\otimes C[-1]$ is the unique
morphism of complexes of degree -1 such that
$\tilde{\Delta}\circ s^{-1} = (s^{-1}\otimes s^{-1})\Delta$.  Here
$s:C\longrightarrow C[-1]$ is the desuspension.

We have $d^2=b^{\prime^2} = db^\prime +b^\prime d=0$, so that
$d+b^\prime$ is a differential of degree -1.  It is moreover a graded
derivation of $B^cC$.

Let $A_1$ and $A_2$ be unital algebras (no differential, no grading).
The {\em shuffle map} is a morphism of $DG$ coalgebras
\begin{equation*}
S:BA_1\otimes BA_2\longrightarrow B(A_1\otimes A_2)\tag{2}
\end{equation*}
defined as follows.  Let $\theta_i:BA_i\to A_i$, $i=1,2$, be universal
twisting cochains of $A_i$.  Define a twisting cochain
\[ \theta :BA_1\otimes BA_2\to A_1\otimes A_2 \]
by
\[ \theta = \theta_1\otimes\varepsilon\eta
+\varepsilon\eta\otimes\theta_2\; , \]
where $\eta$ is the counit map of $BA_i$ and $\varepsilon$ is the unit
map of $A_i$.  Note that $\theta$ is zero except on linear span of tensors
of the type $a_1\otimes 1$, or $1\otimes a_2$.  Checking the twisting
cochain condition for $\theta$ is easy.  By universal property of the bar
construction we obtain $S$.

A simple computation using (1) gives the following explicit formula for
$S$:
\[ S((a_1,\dots , a_p)\otimes (b_1,\dots b_q)) = \sum_{\sigma\in
S_{p,q}}sgn(\sigma )\sigma (a^\prime_1,\dots ,a^\prime_p,b^\prime_1,\dots
,b^\prime_q)\; , \]
where $a^\prime_i=a_i\otimes 1$ and $b^\prime_i=1\otimes b_i$, and
$S_{p,q}$ is the set of all $(p,q)$ shuffle permutations in the symmetric group
$S_{p+q}$.

We need to know that $S$ is a quasi-isomrophism.  This is trivial for
unital algebras over a field, since it is well known that $b^\prime$
is acyclic for unital algebras and one can use the K\"{u}nneth formula
for tensor product of complexes over a field.  We mention, however, that
$S$ has an explicit homotopy inverse
\[ A:B(A_1\otimes A_2)\to BA_1\otimes BA_2\; , \]
the {\em Alexander-Whitney map} (see, e.g. [M]), which proves that
$S$ is a quasi-isomorphism over any ground ring.  Same proof works when
$A_1$ and $A_2$ are unital $DG$ algebras.  It is important, however, to
note that $A$ is not a morphism of $DG$ coalgebras.  This makes finding
an explicit formula for the inverse of the product map in periodic cyclic
homology a difficult task.  See however [B] where this problem is
successfully solved.  An alternative approach would be to extend $A$ to
an $A_\infty$-morphism of $DG$ coalgebras which in theory one knows to
exist and let it act on cyclic complexes.

\section{Periodic cyclic homology of $DG$ coalgebras}
Periodic cyclic homology of $DG$ algebras was defined by Goodwillie and
others (see [G] and references therein).  One of the main results of [G]
is the fact that if $f:A\to B$ is a morphism of $DG$ algebras which is a
quasi-isomorphism of complexes, then the induced map on periodic cyclic
homology $f_\ast :HP_\ast (A)\to HP_\ast (B)$ is an isomorphism.  We need
the analogue of this result for $DG$ coalgebras.  We also need the fact that
for $DG$ coalgebras of finite cohomological dimension (1 and 2 in our
applications), the periodic cyclic complex is quasi-isomorphic to some
higher order versions of the $X$-complex.

Let $C$ be a $DG$ coalgebra.  In [Kh], we have carefully defined the
$DG$ coalgebra of (noncommutative) differential forms over $C$, and
denoted it by $(\Omega C,d)$.  Its construction is dual to the
corresponding construction for algebras and $DG$ algebras.

We need to adopt some basic definitions and constructions from [CQ$_1$] to
our $DG$ coalgebraic set up.  Let $C$ be a $DG$ coalgebra and let
$(\Omega C,d)$ denote the $DG$ coalgebra of {\em universal
codifferential forms} over $C$.  Let $\eta :C\longrightarrow k$ be the
counit of $C$.  We have $\Omega^nC = C\otimes \overline{C}^{\otimes n}$,
where $\overline{C} = Ker\;\eta$.  Let $b:\Omega^\bullet
C\longrightarrow\Omega^{\bullet +1}C$ be the analogue of the {\em
Hochschild boundary operator} and let $N$ be the {\em number operator}
which multiplies a differential form by its degree.  Let
\[ \Omega^{norm}C = ker\{ (b+dN)^2:\Omega C\longrightarrow\Omega C\}\; .
\]

Equipped with the differential $b+dN$ and with its natural
$\mathbb{Z}/2$ grading, $(\Omega^{norm}C, b+dN)$ can be regarded as a
supercomplex.  There is a decreasing filteration
$\{F^n\Omega^{norm}C\}_{n\geq 2}$ on $\Omega^{norm}C$, where $F^n$
consists of forms of degree at least $n$.  The successive quotient
complexes $\Omega^{norm}C/F^n$ approximate the normalized cyclic
bicomplex for $DG$ coalgebras.  We need only the first two quotients,
denoted by $X(C)$ and $X^2(C)$.  These are the supercomplexes
\[ \begin{array}{rl}
X(C): & \begin{CD} C @<\overset{b}{\longrightarrow}<{d}<
\Omega^1C_\natural\end{CD}\\
X^2(C): & \begin{CD} C\bigoplus\Omega^2C_\natural
@<\overset{b+2d}{\longrightarrow}<{b+d}< \dot{\Omega}^1C,
\end{CD}\end{array} \]
where $\natural$ denotes the cocommutator subspace and $\dot{\Omega}^C =
\Omega^{norm,1}C$.  Note that
$\Omega^1C_\natural\subset\dot{\Omega}^1C$.

We are mostly interested in the total complexes of these bicomplexes which
we denote by $\widehat{X}(C)$ and $\widehat{X}^2(C)$.  Here total means
taking direct products.  When working with $X$-complex of $DG$ algebras,
total means taking direct sums.

We define the {\em periodic cyclic complex} of a $DG$ coalgebra $C$ by
\[ CC_\ast (C) = \widehat{X}(B^cC)\; . \]
We denote the homology of this complex by
$HP_\ast (C)$.  The motivation for this definition is as follows.  In
[Q$_1$], D. Quillen computed $\widehat{X}(BA)$, when $A$ is an algebra, and
showed that it is isomorphic with the total complex of Connes-Tsygan
bicomplex of $A$.  The same proof easily extends to the case of $DG$
algebras.  Our definition of periodic cyclic complex of $DG$ coalgebras
is simply the dual of this definition-theorem for $DG$ algebras.

The following result was first proved by T. Goodwillie for $DG$ algebras
[G].

\begin{proposition}
Let $f:C\to D$ be a morphism of $DG$ coalgebras such that $f$ is a
quasi-isomorphism of complexes.  Then $f_\ast :HP_\ast (C)\to HP_\ast
(D)$ is an isomorphism.
\end{proposition}

Let $C$ be a $DG$ coalgebra.  Note that we have a natural morphism of
complexes
\[ I:\widehat{X}(C)\longrightarrow\widehat{X}^2(C)\; , \]
obtained from the inclusion $C\longrightarrow C\oplus\Omega^2C_\natural$
and $\Omega^1C_\natural\longrightarrow\dot{\Omega}^1C$.  In general,
there is no natural map $\widehat{X}^2(C)\longrightarrow\widehat{X}(C)$.
However, it is shown in [Kh] that if $C = BA$ is the bar construction
then $I$ is a homotopy equivalence and a homotopy inverse
\[ R:\widehat{X}^2(BA)\longrightarrow\widehat{X}(BA) \]
is constructed.  There is a dual statement for $DG$ coalgebras, but we
don't need it in this paper.

Let us call a $DG$ coalgebra of {\em finite cohomological dimension} if
its underlying coalgebra has finite cohomological dimension.  We need to
know that if a $DG$ coalgebra $C$ has finite cohomological
dimension, then its periodic cyclic complex $CC_\ast (C)$ is
quasi-isomorphic to a ``small'' complex.  We need this only for
cohomological dimensions 1 and 2.

\begin{proposition}
Let $C$ be a $DG$ coalgebra.  Then
\begin{enumerate}
\item[1.] The natural map $\widehat{X}(C)\longrightarrow CC_\ast (C)$
is a quasi-isomorphism if $C$ has cohomolgical dimension 1.
\item[2] The natural map $\widehat{X}^2(C)\longrightarrow CC_\ast (C)$
is a quasi-isomorphism if $C$ has cohomological dimension 2.
\end{enumerate}
\end{proposition}

\section{A Proof of the K\"{u}nneth formula}
Let $A_1$ and $A_2$ be unital algebras over a field of characteristic
zero.  In general, $X(A_1)\otimes X(A_2)$ can not be quasi-isomorphic to
$X(A_1\otimes A_2)$.  However, when $A_1$ and $A_2$ are quasifree (in
particular free) algebras then the cohomological dimension of $A_1\otimes
A_2$ is at most 2 and $X(A_1)\otimes X(A_2)$ is quasi-isomorphic to
$X^2(A_1\otimes A_2)$.  This fact is due to Cuntz and Quillen [CQ$_1$].  In
order to obtain explicit formulas, M. Puschnigg constructed a morphism
of complexes [P],
\[ P:X^2(A_1\otimes A_2)\longrightarrow X(A_1)\otimes X(A_2)\; , \]
and showed that $P$ is a quasi-isomorphism when $A_1$ and $A_2$ are
quasi-free.  In fact, an explicit right inverse to $P$ is the map defined in
[CQ$_1$], for arbitrary algebras, $X(A_1)\otimes X(A_2)\longrightarrow
X(A_1\otimes A_2)$ combined with the canonical inclusion
$X(A_1\otimes A_2)\longrightarrow X^2(A_1\otimes A_2)$.

Next we observe that the above constructions are completely functorial
and extend to $DG$ algebras and $DG$ coalgebras.  To be precise, let us
define the {\em completed tensor product} $C\widehat{\otimes} D$ of
infinite product vector spaces \d C = \prod_{i\geq 0} C_i$ and \d D =
\prod_{i\geq 0}D_i$ by
\[
C\widehat{\otimes}D = \prod_{n\geq 0}\bigoplus_{i_+j=n}C_i\otimes D_j\; .
\]
There is an obvious injection
$C\otimes D\longrightarrow C\widehat{\otimes}D$.  Since the bar construction
is free as a coalgebra, dualizing the above map $P$ we obtain a
quasi-isomorphism of complexes
\[ \widehat{X}(BA_1)\widehat{\otimes}\widehat{X}(BA_2)\longrightarrow
\widehat{X}^2(BA_1\otimes BA_2)\; . \]

Consider the sequence of maps
\[ \widehat{X}(BA_1)\widehat{\otimes}\widehat{X}(BA_2)\overset{P}{\longrightarrow}
\widehat{X}^2(BA_1\otimes BA_2)\overset{\overline{S}}{\longrightarrow}
\widehat{X}^2(B(A_1\otimes A_2))\overset{R}{\longrightarrow}
\widehat{X}(B(A_1\otimes A_2))\; , \]
where $\overline{S}$ is induced by the shuffle map
$S:BA_1\otimes BA_2\longrightarrow B(A_1\otimes A_2)$ and $R$
is the retraction introduced earlier.
\begin{theorem}
$P$, $\overline{S}$ and $R$ are quasi-isomorphisms.
\end{theorem}
\begin{proof}
We only have to show that $\overline{S}$ is a quasi-isomorphism.
Consider the diagram
\[ \begin{CD}
\widehat{X}^2(BA_1\otimes BA_2) @<<< \widehat{X}(B^c(BA_1\otimes BA_2) =
CC_\ast (BA_1\otimes BA_2)\\
@VV \overline{S} V @ VV\widetilde{S} V\\
\widehat{X}^2(B(A_1\otimes A_2)) @<<< \widehat{X}(B^cB(A_1\otimes
A_2))=CC_\ast (B(A_1\otimes A_2))
\end{CD} \]
where $\widetilde{S}$ is the map induced from the $DG$ coalgebra map
$S:BA_1\otimes BA_2\longrightarrow B(A_1\otimes A_2)$ on periodic
cyclic complexes.  By proposition 2.1, $\widetilde{S}$ is a
quasi-isomorphism as $S$ is a quasi-isomorphism of complexes.  By
proposition 2.2, the horizontal arrows are also quasi-isomorphisms and hence
$\overline{S}$ is a quasi-isomorphism.
\end{proof}

It follows that the composition
\[ \widehat{S} = R\circ \overline{S}\circ P:CC_\ast
(A_1)\widehat{\otimes}CC_\ast (A_2)\longrightarrow CC_\ast (A_1\otimes A_2) \]
is a quasi-isomorphism of complexes.

\vskip30pt

\noindent Masoud Khalkhali\\
University of Western Ontario\\
London, Canada\\
N6A 5B7\\
masoud@julian.uwo.ca

\end{document}